\documentclass[11pt]{amsart}
\usepackage[margin=1.1in]{geometry}
\usepackage{amsmath,amssymb,amsthm}
\usepackage{hyperref}

\newtheorem{theorem}{Theorem}[section]
\newtheorem{lemma}[theorem]{Lemma}
\newtheorem{proposition}[theorem]{Proposition}
\theoremstyle{remark}
\newtheorem{remark}[theorem]{Remark}
\renewcommand{\sec}{\operatorname{sec}}

\title{An isolated totally geodesic flat torus in a metric of
nonnegative curvature on $S^2\times S^2$}

\begin{document}
\author{Dmitri Burago}

\maketitle
\begin{abstract}
Let $E$ be a triaxial ellipsoid with its induced metric, and let
$\gamma$ be the principal ellipse perpendicular to the middle axis.
In the product $E\times E$, the torus $\gamma\times\gamma$ is flat
and totally geodesic. We prove that some neighbourhood of this torus
contains no other embedded totally geodesic flat torus. We also show
that the linearized totally geodesic equation has trivial normal kernel.
The proof uses the hyperbolicity of $\gamma$ and the splitting of flat
totally geodesic tori in products of positively curved surfaces.
\end{abstract}

\section{Motivation}

Bruce Kleiner raised this question in discussions of the H. Hopf conjecture.
The diagonal Cheeger deformation of the round product reduces the set
of flat planes and retains families of totally geodesic flat tori
\cite[\S\S2 and 4.4]{Bettiol}. Is it possible for such a torus to
be isolated? For recent work on positive sectional curvature on
$S^2\times S^2$, see \cite[Theorem~1.1]{BrendleHung}. This elementary note arose from a separate investigation of curvature
variations involving jet computations assisted by Claude Max. Verified by Astra 6 pro. 

In dimension three, a two-sided totally geodesic flat torus in a
nonnegatively curved manifold has flat equidistants. The following
observation explains what persists under normal translation.

\begin{proposition}\label{prop:equidistant}
Let $T$ be a two-sided totally geodesic flat torus in a Riemannian
three-manifold with $\sec\geq0$. For sufficiently small $|r|$, the
equidistant torus $T_r$ is intrinsically flat, its tangent planes
have ambient sectional curvature zero, and its shape operator is
semidefinite of rank at most one.
\end{proposition}

\begin{proof}
Use the normal exponential map to define $T_r$. With the convention
$S_rX=\nabla_X\partial_r$, the shape operator satisfies the Riccati
equation
\[
 S_r'+S_r^2+R_{\partial_r}=0,\qquad
 R_{\partial_r}X=R(X,\partial_r)\partial_r,\qquad S_0=0.
\]
Use parallel transport along the normal geodesics to compare the
operators. Since $R_{\partial_r}\geq0$, we have $S_r'\leq0$. Hence
$S_r\leq0$ for $r>0$ and $S_r\geq0$ for $r<0$. In either case
$\det S_r\geq0$. The Gauss equation gives
\[
 K_{T_r}=\sec(TT_r)+\det S_r.
\]
Both terms on the right are nonnegative. Their sum has integral zero
by Gauss--Bonnet, so both vanish pointwise. The determinant condition
and semidefiniteness give the claimed rank bound.
\end{proof}

The shape operator of an equidistant torus can have rank one.
The question here concerns isolation among tori whose entire second
fundamental form vanishes.

\section{The example}

Let $a>b>c>0$ and set
\[
 E=\left\{(x,y,z)\in\mathbb R^3:
       \frac{x^2}{a^2}+\frac{y^2}{b^2}+\frac{z^2}{c^2}=1
       \right\},\qquad
 \gamma=E\cap\{y=0\}.
\]
The reflection $(x,y,z)\mapsto(x,-y,z)$ is an isometry of $E$.
Its fixed curve $\gamma$ is therefore a simple closed geodesic.
Write $M=E\times E$ and $T=\gamma\times\gamma$.

For a totally geodesic surface $T\subset M$, let $\mathcal J(T)$
be the space of smooth normal fields $V$ satisfying
\begin{equation}\label{eq:linearized}
 \nabla_X^\perp\nabla_Y^\perp V
 -\nabla_{\nabla_X^T Y}^\perp V
 +(R(V,X)Y)^\perp=0
 \qquad(X,Y\in TT).
\end{equation}
Here $\nabla^\perp$ is the normal connection, and our curvature
convention is $\langle R(X,Y)Y,X\rangle=\sec(X,Y)|X\wedge Y|^2$.
Equation \eqref{eq:linearized} is the linearization of the vanishing
of the second fundamental form under a normal variation. In
particular, a variation through totally geodesic surfaces has
normal velocity in $\mathcal J(T)$.

\begin{theorem}\label{thm:main}
The product metric on $M\cong S^2\times S^2$ has nonnegative
sectional curvature, and $T$ is a totally geodesic flat torus.
Moreover, some open neighbourhood of $T$ contains no other
embedded totally geodesic flat torus, and $\mathcal J(T)=0$.
\end{theorem}

\section{Splitting in a product}

If $\Sigma_1$ and $\Sigma_2$ have positive Gauss curvature, then
for $X=(X_1,X_2)$ and $Y=(Y_1,Y_2)$,
\begin{equation}\label{eq:sec}
 \langle R(X,Y)Y,X\rangle
 =K_1|X_1\wedge Y_1|^2+K_2|X_2\wedge Y_2|^2.
\end{equation}
Hence their product has nonnegative curvature. A plane has zero
curvature precisely when both projections have rank one. It is
then the direct sum of a line in the first factor and a line in
the second. We call such a plane mixed. Since $E$ has positive
Gauss curvature, this proves the curvature assertion of the
theorem. The product of two geodesics is totally geodesic and has
the flat product metric, so it also proves the assertion about $T$.

\begin{lemma}\label{lem:class}
Every compact connected embedded totally geodesic flat torus
$T'\subset\Sigma_1\times\Sigma_2$ is $c\times d$ for embedded
closed geodesics $c\subset\Sigma_1$ and $d\subset\Sigma_2$.
\end{lemma}

\begin{proof}
The Gauss equation and \eqref{eq:sec} show that every tangent plane
of $T'$ is mixed. Thus each factor projection has rank one.
Its fibres are compact one-dimensional submanifolds by the
constant rank theorem. A component of a fibre is a circle tangent
to one factor. It is a geodesic in that factor. Indeed, if $X$ is
a local unit tangent along the fibre, then $\nabla_X X$ lies
both in $TT'$ and in the tangent bundle to that factor. Their
intersection is the line spanned by $X$, and $\nabla_X X$ is
orthogonal to $X$. Embeddedness of $T'$ makes the fibre component
an embedded circle in its factor.

Through a point $(p,q)$ of $T'$, choose the two fibre circles
$c\times\{q\}$ and $\{p\}\times d$. The totally geodesic surfaces
$T'$ and $c\times d$ have the same tangent plane there, so their
germs agree by the exponential map. The set of points in
$c\times d$ where the germs agree is open. It is also closed:
at a limit of such points, compactness puts the point in $T'$,
and continuity gives equality of the tangent planes. The germ
argument then applies again. Since $c\times d$ is connected,
it is contained in $T'$. It is open and closed there, so
connectedness of $T'$ gives $T'=c\times d$.
\end{proof}

\section{Isolation}

\begin{lemma}\label{lem:hyp}
The closed geodesic $\gamma$ is hyperbolic. In particular, if $K(s)$
is the Gauss curvature along its arclength parametrization, the
equation
\[
 u''+K(s)u=0
\]
has no nonzero solution periodic with the period of $\gamma$.
\end{lemma}

\begin{proof}
The hyperbolicity of the middle principal ellipse of a triaxial
ellipsoid is a classical result about the ellipsoid geodesic flow.
See \cite[Theorem~3.5.16]{KliRG} and the discussion in
\cite[\S1.5.2]{CG}. The closed geodesic through the umbilics is
the ellipse with semiaxes $a$ and $c$, as described by Arnold
\cite[\S47]{Arnold}. The monodromy of the
scalar Jacobi equation is the linearized Poincar\'e return map
in the normal directions. Hyperbolicity excludes the eigenvalue
$1$, giving the second assertion.
\end{proof}

\begin{lemma}\label{lem:tube}
For every neighbourhood $N$ of the two oriented unit tangent
lifts of $\gamma$, there is a tubular neighbourhood $V$ of
$\gamma$ such that every unit-speed complete geodesic contained
in $V$ has its entire unit tangent lift in $N$.
\end{lemma}

\begin{proof}
Otherwise choose complete geodesics in tubes of radii tending to
zero and unit tangent vectors along them outside $N$. A convergent
subsequence of these vectors gives an initial vector outside $N$.
By continuity of the geodesic flow, the resulting complete geodesic
lies in $\gamma$ at every time. Its initial vector is therefore
one of the two unit vectors tangent to $\gamma$ at its base point,
a contradiction.
\end{proof}

\begin{proof}[Proof of positional isolation in Theorem~\ref{thm:main}]
Each oriented lift of $\gamma$ is a hyperbolic periodic orbit.
It has an isolating neighbourhood in the unit tangent bundle:
any orbit that stays in the neighbourhood for all time is the
periodic orbit itself. This follows by taking a transverse section
and using the local stable and unstable manifolds of the hyperbolic
fixed point of the return map \cite[\S3.3]{KliRG}. Choose disjoint isolating
neighbourhoods of the two orientations and apply Lemma~\ref{lem:tube}
to their union.

Every closed geodesic contained in the resulting tube $V$ has
image $\gamma$. A totally geodesic flat torus in $V\times V$ is
$c\times d$ by Lemma~\ref{lem:class}. Both $c$ and $d$ are contained
in $V$, and hence $c=d=\gamma$ as subsets.
\end{proof}

\begin{proof}[Proof of infinitesimal isolation in Theorem~\ref{thm:main}]
Use arclength coordinates $(s_1,s_2)$ and the parallel normal frame
$(\nu_1,\nu_2)$ induced by the unit normals to $\gamma$ in each
factor. Write $V=u\nu_1+w\nu_2$. The product connection and product
curvature turn \eqref{eq:linearized} into
\begin{align*}
 \partial_1^2u+K(s_1)u&=0, & \partial_1^2w&=0,\\
 \partial_2^2u&=0, & \partial_2^2w+K(s_2)w&=0,\\
 \partial_1\partial_2u&=0, & \partial_1\partial_2w&=0.
\end{align*}
Periodicity implies that $u$ is independent of $s_2$ and $w$ is
independent of $s_1$. The two remaining equations are the scalar
Jacobi equation along $\gamma$, with periodic boundary conditions.
Lemma~\ref{lem:hyp} gives $u=w=0$.
\end{proof}

\section{Scope of the example}

\begin{remark}[Embedded and immersed geodesics]
For isolation among embedded closed geodesics, nondegeneracy
already suffices. Let $\gamma$ be a simple nondegenerate closed
geodesic in an orientable surface. Lemma~\ref{lem:tube} applies
to a sufficiently thin annular tube about $\gamma$. In that tube,
every embedded closed geodesic projects locally diffeomorphically
onto $\gamma$ under normal projection. An essential embedded circle
in an annulus has winding number $1$ or $-1$, so such a geodesic
makes exactly one circuit. Its unit tangent at a fixed transverse
fibre therefore gives a fixed point of the one-circuit Poincar\'e
map, near the fixed point for $\gamma$. Nondegeneracy makes that
fixed point isolated. By Lemma~\ref{lem:class}, the same argument
applies to the product torus formed from two such geodesics.

Closed geodesics that wind several times around an elliptic orbit
may arise in a thin tube. Such curves cannot be embedded in that
annulus and therefore do not establish accumulation of embedded
flat tori. The hyperbolicity argument above gives the stronger
confinement statement for all complete geodesics in the tube.
\end{remark}

\begin{remark}[Zero-curvature planes]
Every mixed plane in $E\times E$ has zero curvature. A stronger
question is whether one can find a nonnegatively curved metric on
$S^2\times S^2$ whose zero-curvature planes are exactly the tangent
planes of one totally geodesic flat torus. One can ask this with
the additional requirement that the metric be close to the round
product. The example above answers the isolation question for
embedded tori. It leaves these stronger questions open.
\end{remark}

\end{document}